\documentclass[12pt, reqno]{amsart}
\usepackage{amsmath}
\usepackage{amsfonts}
\usepackage[all]{xy}
\begin{document}
\numberwithin{equation}{section}

\def\1#1{\overline{#1}}
\def\2#1{\widetilde{#1}}
\def\3#1{\widehat{#1}}
\def\4#1{\mathbb{#1}}
\def\5#1{\frak{#1}}
\def\6#1{{\mathcal{#1}}}

\newcommand{\de}{\partial}
\newcommand{\R}{\mathbb R}
\newcommand{\Ha}{\mathbb H}
\newcommand{\al}{\alpha}
\newcommand{\tr}{\widetilde{\rho}}
\newcommand{\tz}{\widetilde{\zeta}}
\newcommand{\tk}{\widetilde{C}}
\newcommand{\tv}{\widetilde{\varphi}}
\newcommand{\hv}{\hat{\varphi}}
\newcommand{\tu}{\tilde{u}}
\newcommand{\tF}{\tilde{F}}
\newcommand{\debar}{\overline{\de}}
\newcommand{\Z}{\mathbb Z}
\newcommand{\C}{\mathbb C}
\newcommand{\Po}{\mathbb P}
\newcommand{\zbar}{\overline{z}}
\newcommand{\G}{\mathcal{G}}
\newcommand{\So}{\mathcal{S}}
\newcommand{\Ko}{\mathcal{K}}
\newcommand{\U}{\mathcal{U}}
\newcommand{\B}{\mathbb B}
\newcommand{\oB}{\overline{\mathbb B}}
\newcommand{\Cur}{\mathcal D}
\newcommand{\Dis}{\mathcal Dis}
\newcommand{\Levi}{\mathcal L}
\newcommand{\SP}{\mathcal SP}
\newcommand{\Sp}{\mathcal Q}
\newcommand{\A}{\mathcal O^{k+\alpha}(\overline{\mathbb D},\C^n)}
\newcommand{\CA}{\mathcal C^{k+\alpha}(\de{\mathbb D},\C^n)}
\newcommand{\Ma}{\mathcal M}
\newcommand{\Ac}{\mathcal O^{k+\alpha}(\overline{\mathbb D},\C^{n}\times\C^{n-1})}
\newcommand{\Acc}{\mathcal O^{k-1+\alpha}(\overline{\mathbb D},\C)}
\newcommand{\Acr}{\mathcal O^{k+\alpha}(\overline{\mathbb D},\R^{n})}
\newcommand{\Co}{\mathcal C}
\newcommand{\Hol}{{\sf Hol}(\mathbb H, \mathbb C)}
\newcommand{\Aut}{{\sf Aut}(\mathbb D)}
\newcommand{\D}{\mathbb D}
\newcommand{\oD}{\overline{\mathbb D}}
\newcommand{\oX}{\overline{X}}
\newcommand{\loc}{L^1_{\rm{loc}}}
\newcommand{\la}{\langle}
\newcommand{\ra}{\rangle}
\newcommand{\thh}{\tilde{h}}
\newcommand{\N}{\mathbb N}
\newcommand{\kd}{\kappa_D}
\newcommand{\Hr}{\mathbb H}
\newcommand{\ps}{{\sf Psh}}
\newcommand{\Hess}{{\sf Hess}}
\newcommand{\subh}{{\sf subh}}
\newcommand{\harm}{{\sf harm}}
\newcommand{\ph}{{\sf Ph}}
\newcommand{\tl}{\tilde{\lambda}}
\newcommand{\gdot}{\stackrel{\cdot}{g}}
\newcommand{\gddot}{\stackrel{\cdot\cdot}{g}}
\newcommand{\fdot}{\stackrel{\cdot}{f}}
\newcommand{\fddot}{\stackrel{\cdot\cdot}{f}}
\def\v{\varphi}
\def\Re{{\sf Re}\,}
\def\Im{{\sf Im}\,}
\def\ext{{\sf ext}\,}
\def\Lr{{\sf Lr}\,}
\def\tr{{\sf tr}\,}

\def\Label#1{\label{#1}}


\def\cn{{\C^n}}
\def\cnn{{\C^{n'}}}
\def\ocn{\2{\C^n}}
\def\ocnn{\2{\C^{n'}}}
\def\je{{\6J}}
\def\jep{{\6J}_{p,p'}}
\def\th{\tilde{h}}


\def\dist{{\rm dist}}
\def\const{{\rm const}}
\def\rk{{\rm rank\,}}
\def\id{{\sf id}}
\def\aut{{\sf aut}}
\def\Aut{{\sf Aut}}
\def\CR{{\rm CR}}
\def\GL{{\sf GL}}
\def\Re{{\sf Re}\,}
\def\Im{{\sf Im}\,}
\def\U{{\sf U}}
\def\Tang{{\rm Tang}_1(\mathbb{C}^N,0)}
\def\la{\langle}
\def\ra{\rangle}

\emergencystretch15pt \frenchspacing

\newtheorem{theorem}{Theorem}[section]
\newtheorem*{theorem**}{Theorem \mynumber}
\newenvironment{theorem*}[1]
  {\newcommand{\mynumber}{\ref{#1}}\begin{theorem**}}
  {\end{theorem**}}
\newtheorem{lemma}[theorem]{Lemma}
\newtheorem{proposition}[theorem]{Proposition}
\newtheorem{problem}[theorem]{Problem}
\newtheorem{corollary}[theorem]{Corollary}

\theoremstyle{definition}
\newtheorem{definition}[theorem]{Definition}
\newtheorem{example}[theorem]{Example}

\theoremstyle{remark}
\newtheorem{remark}[theorem]{Remark}
\numberwithin{equation}{section}

\title[Resonances in Loewner equations]{Resonances in Loewner equations}
\author[L. Arosio]{Leandro Arosio}
\address{Dipartimento di Matematica, Universit\`{a} di Roma \textquotedblleft La Sapienza\textquotedblright, Piazzale A. Moro 5, 00185 Roma,
Italy. }
\email{arosio@mat.uniroma1.it}

\date{\today }


\setcounter{tocdepth}{1}

\begin{abstract}
We prove that given a  Herglotz vector field on the unit ball of $\mathbb{C}^n$ of the form $H(z,t)=(a_1 z_1,\dots,a_n z_n)+O(|z|^2)$ with $\Re a_j<0$ for all $j$, its evolution family admits an associated Loewner chain, which is normal if no real resonances occur. Hence the Loewner-Kufarev PDE  admits a solution defined for all positive times.
\end{abstract}
\maketitle
\tableofcontents

\section{Introduction}

Classical Loewner theory in the unit disc $\mathbb{D}\subset \mathbb{C}$ was introduced by C. Loewner in 1923 \cite{Loewner} and developed with contributions of P.P. Kufarev in 1943 \cite{Kufarev} and C. Pommerenke in 1965 \cite{Pommerenke}, and has been since then used to prove several deep results in geometric function theory \cite{Graham-Kohr}.  Loewner theory is  one of the main ingredients of the proof of the  Bieberbach conjecture given by de Branges \cite{deBranges} (see also \cite{Fitzgerald-Pommerenke}) in 1985.

Among the extensions of classical Loewner theory we recall the chordal Loewner theory \cite{Kufarev-Sobolev-Sporysheva}, the celebrated theory of Schramm-Loewner evolution \cite{Schramm} introduced in 1999 and the theory of Loewner chains in several complex variables  \cite{Pfaltzgraff}\cite{Duren-Graham-Hamada-Kohr}\cite{Graham-Hamada-Kohr-Kohr}. 

In  \cite{Bracci-Contreras-Diaz}\cite{Bracci-Contreras-Diaz-II} it is proposed  a generalization of both the radial and chordal theories. It is shown that on complete hyperbolic manifolds there is a one-to-one correspondence between certain semicomplete non-autonomous holomorphic vector fields (called {\sl Herglotz vector fields} and denoted $H(z,t)$) and families $(\v_{s,t})_{0\leq s\leq t}$ of holomorphic self-maps called {\sl evolution families}.
Indeed,  if  $H(z,t)$ is a Herglotz vector field, then the family $(\v_{s,t})$ of evolution operators for the Loewner-Kufarev ODE 
\begin{equation}
\overset{\bullet}z(t)=H(z,t),\quad t\geq 0,\  z\in\mathbb{B},
\end{equation}
 is an evolution family. Conversely, any  evolution family is the family of evolution operators for some Loewner-Kufarev ODE. 

In \cite{Contreras-Diaz-Gumenyuk} it is proved that in  dimension one evolution families are (up to biholomorphism) in one-to-one correspondence with image-growing families  $(f_s)_{s \geq 0}$ of univalent mappings $f_s\colon\mathbb{D}\rightarrow \mathbb{C}$ called {\sl Loewner chains}. Namely given any Loewner chain $(f_s)$ the family $(\v_{s,t})$ defined by $$\v_{s,t}=f_t^{-1}\circ f_s$$ is an evolution family, which is said to be {\sl associated} to $(f_s)$. Conversely given any evolution family, there exists an associated Loewner chain.
Composing the two correspondences above we obtain the correspondence between  Loewner chains and Herglotz vector fields: $H(z,t)$ is associated to $(f_s)$ if and only if the mapping $t\mapsto f_t$ is a global solution for the  Loewner-Kufarev PDE
\begin{equation}
\frac{\partial f_t(z)}{\partial t}=-f'_t(z)H(z,t),\quad t\geq 0,\  z\in\mathbb{D}.
\end{equation}

Let $N$  be an integer greater or equal to $2$, and let $\mathbb{B}$ be the unit ball of $\mathbb{C}^N$.  A Loewner chain on $\mathbb{B}$ is an image-growing family  $(f_s)_{s \geq 0}$ of univalent mappings $f_s\colon\mathbb{\mathbb{B}}\rightarrow \mathbb{C}^N$. Every Loewner chain admits an associated evolution family, but it is not known whether the converse is true.
In \cite{Arosio-Bracci-Hamada-Kohr} it is proposed an abstract approach to the notion of Loewner chain. Let $M$ be an $N$-dimensional complete hyperbolic complex manifold. An abstract Loewner chain is an image-growing family $(f_s)$ of univalent mappings defined on $M$ which are allowed to take values on an arbitrary $N$-dimensional complex manifold. In \cite{Arosio-Bracci-Hamada-Kohr} it is shown that to any evolution family $(\v_{s,t})$ on $M$ there corresponds a unique (up to biholomorphisms) abstract Loewner chain $(f_s)$.  In this way one can define the {\sl Loewner range manifold} of $(\v_{s,t})$ $$\Lr(\v_{s,t})=\bigcup_{s\geq 0}f_s(M)$$ which is well defined and unique up to biholomorphism. Hence the classical problem of finding a Loewner chain (with values in $\mathbb{C}^N$) associated to a given evolution family $(\v_{s,t})$ of the unit ball $\mathbb{B}\subset\mathbb{C}^N$ corresponds to investigating whether the Loewner range manifold $\Lr(\v_{s,t})$ embeds holomorphically in $\mathbb{C}^N$. 

In this paper we investigate this problem for a special type of evolution families  on  $\mathbb{B}$.
Let $\Lambda $ be an  ($N\times N$)-complex matrix
\begin{equation}\label{dilationmatrix}
\Lambda =\mbox{Diag}(\alpha_1,\dots,\alpha_N),\quad \mbox{where}\quad\Re{\alpha_N}\leq\dots\leq\Re{\alpha_1}<0.
\end{equation}
We define a {\sl dilation evolution family} as an evolution family $(\v_{s,t})$ on the unit ball $\mathbb{B}\subset \mathbb{C}^N$ satisfying $$ \varphi_{s,t}(z)=e^{\Lambda ( t-s)}z+O(|z|^2).$$
A {\sl normal Loewner chain} is a Loewner chain $(f_s)$ such that
\begin{enumerate}
\item $f_s(z)=e^{-\Lambda s}z+O(|z|^2)$,
\item $(h_s)$  is a normal family.
\end{enumerate}
Notice that  each $h_s=e^{\Lambda s} f_s$ fixes the origin and is tangent to identity in the origin (we say $h_s\in \Tang$).

\begin{problem}\label{question}
Given a  dilation evolution family, does there exist an associated Loewner chain (with values in $\mathbb{C}^N$)?
\end{problem}
An affirmative answer to Problem \ref{question} would yield as a consequence that any Loewner-Kufarev  PDE 
\begin{equation}
\frac{\partial f_t(z)}{\partial t}=-d_zf_tH(z,t),\quad t\geq 0,\  z\in\mathbb{B},
\end{equation}
where $H(z,t)=\Lambda z+ O(|z|^2)$ (in this case the equation is known as the Loewner PDE), admits  global solutions. A partial answer may be obtained by simply combining   \cite[Theorem 3.1]{Duren-Graham-Hamada-Kohr} and   \cite[Theorems 2.3, 2.6]{Graham-Hamada-Kohr-Kohr}:
\begin{theorem}\label{kohrresult}
Let $(\varphi_{s,t})$ be a dilation evolution family such that the eigenvalues of $\Lambda $ satisfy 
\begin{equation}\label{condition}
2\Re\alpha_1<\Re\alpha_N.
\end{equation}
Then there exists a normal Loewner chain $(f_s)$ associated to $(\varphi_{s,t})$, such that $\bigcup_s f_s(\mathbb{B})=\mathbb{C}^N,$ hence  $\Lr(\v_{s,t})=\mathbb{C}^N$. This chain is given by 
\begin{equation}\label{koenigsiteration}
f_s=\lim_{t\rightarrow +\infty}e^{-\Lambda t}\v_{s,t},
\end{equation}
where the limit is taken in the topology of uniform convergence on compacta,
and it is the unique normal Loewner chain associated to $(\varphi_{s,t})$.
A family of univalent mappings $(g_s)$ is a Loewner chain  associated to $(\varphi_{s,t})$ if and only if there exists an entire univalent mapping $\Psi$ on $\mathbb{C}^{N}$ such that $$g_s=\Psi\circ f_s.$$ 
\end{theorem}

The main result of this paper gives an affirmative answer to Problem \ref{question}, without assuming condition (\ref{condition}).
\begin{theorem*}{continuousresult}
Let $(\varphi_{s,t})$ be a dilation evolution family. Then there exists a Loewner chain $(f_s)$ associated to $(\varphi_{s,t})$, such that $\bigcup_s f_s(\mathbb{B})=\mathbb{C}^N,$ hence  $\Lr(\v_{s,t})=\mathbb{C}^N$. If no real resonances occur among the eigenvalues of $\Lambda$, then $(f_s)$ is a normal chain, not necessarily unique. A family of univalent mappings $(g_s)$ is a Loewner chain  associated to $(\varphi_{s,t})$ if and only if there exists an entire univalent mapping $\Psi$ on $\mathbb{C}^{N}$ such that $$g_s=\Psi\circ f_s.$$ 
\end {theorem*}

Notice that (\ref{condition}) is a classical condition which ensures the existence of a solution for the Schr\"oder functional equation. In fact we will see that normal Loewner chains correspond to solutions of a parametric Schr\"oder equation. Let us first recall some facts about linearization of germs.

Let $\v(z)=e^\Lambda z+O(|z|^2)$ be a holomorphic germ at the origin of $\mathbb{C}^N$, where $\Lambda$ is a matrix satisfying (\ref{dilationmatrix}). If $h\in\Tang$ is a solution of the Schr\"oder equation 
\begin{equation}\label{schroderseveral}
h\circ \v=e^\Lambda h,
\end{equation}
we say that $h$ linearizes $\v$. It is not always possible to solve this equation, indeed there can occur complex resonances among the eigenvalues of $\Lambda$, that is algebraic identities  $$\sum_{j=1}^N k_j\alpha_j=\alpha_l,$$ where $k_j\geq0$ and $\sum_j k_j\geq 2$, which are obstructions to linearization (the term ``complex'' is not standard and is here used to  distinguish from real resonances, defined below). Indeed  a celebrated theorem of Poincar\'e (see for example \cite[pp. 80--86]{Rudin-Rosay})  states that if no complex resonances occur, then there exists a solution $h$ for (\ref{schroderseveral}). If moreover  $2\Re\alpha_1<\Re\alpha_N$ then $h$ is given by $\lim_{n\rightarrow +\infty}e^{-\Lambda n}\v^{\circ n}$. 

In our case we are interested in the following parametric Schr\"oder equation
\begin{equation}\label{lastbutnotleast}
h_m\circ \v_{n,m}=e^{\Lambda (m-n)} h_n,
\end{equation}
where $(\v_{n,m})$ is the discrete analogue of a dilation evolution family. We search for a solution $(h_n)$ which is a  normal family of univalent mappings in $\Tang$. The parametric Schr\"oder equation admits such a solution $(h_n)$ if and only if $(\v_{n,m})$ admits a discrete normal Loewner chain $(f_n)$, and $$(f_n)=(e^{-\Lambda n} h_n).$$

There are surprising differences between the Schr\"oder functional equation (\ref{schroderseveral}) and  (\ref{lastbutnotleast}). Namely, while in the first complex resonances are  obstructions to the existence of formal solutions, in the latter there always exists the holomorphic solution $h_n=e^{\Lambda n} \v_{0,n}^{-1}$, but  the domain of definition of the mapping $h_n$  shrinks as $n$ grows. If, as we need, we look for solutions which are all defined in the unit ball $\mathbb{B}$, then we find as obstructions  {\sl real resonances} among the eigenvalues of $\Lambda$, that is algebraic identities  $$\Re(\sum_{j=1}^N k_j\alpha_j)=\Re\alpha_l,$$ where $k_j\geq0$ and $\sum_j k_j\geq 2$. If real resonances occur we solve a slightly different equation: $$h_m\circ \v_{n,m}=T_{n,m}\circ h_n,$$ where $(T_{n,m})$ is a suitable triangular evolution family, finding this way a non-necessarily normal discrete Loewner chain associated to $(\v_{n,m})$. 

Once we solved the problem for discrete times, we  solve the problem for continuous times: we discretize a given continuous dilation evolution family $(\v_{s,t})$ obtaining a discrete dilation evolution family $(\v_{n,m})$, and we take the associated discrete Loewner chain $(f_n)$. Then we extend $(f_n)$ to all real positive times obtaining this way a Loewner chain $(f_s)$ and Theorem \ref{continuousresult} above.

We give examples of 
\begin{enumerate}
\item a dilation evolution family with no real resonances and several associated normal Loewner chains,
\item a semigroup-type dilation evolution family with complex resonances which does not admit any associated normal Loewner chain,
\item  a discrete dilation evolution family with pure real resonances (real non-complex resonances) which does not admit any discrete normal Loewner chain,
\item a discrete evolution family not of dilation type which does not admit any associated discrete Loewner chain (with values in $\mathbb{C}^N$).
\end{enumerate}

I want to thank Prof. F. Bracci for suggesting the problem and for his precious help. After writing the preliminary version of this paper, I became acquainted with the work of M. Voda ``Solution of a Loewner chain equation in several complex variables'' (arXiv:1006.3286v1 [math.CV], 2010), where an analogue of Theorem \ref{continuousresult} is proved with completely different methods. I want to thank Prof. M. Contreras for informing me about this work. I want to thank M. Voda and the referee for precious comments and remarks.

\section{Preliminaries}\label{preliminaries}
The following is a several variables version of the Schwarz Lemma \cite[Lemma 6.1.28]{Graham-Kohr}.
\begin{lemma}\label{schwarz}
Let $M>0$ and $f\colon \mathbb{B}\rightarrow \mathbb{C}^N$ be a holomorphic mapping fixing the origin and bounded by $M$. Then for $z$ in the ball, $|f(z)|\leq M|z|$. If there is a point $ z_0\in \mathbb{B}\setminus \{0\}$ such that $|f( z_0)|=M| z_0|$, then $|f(\zeta z_0)|=M|\zeta  z_0|$ for all $|\zeta|<1/| z_0|.$ Moreover, if $f(z)=O(|z|^k),\ k\geq 2$, then for $z$ in the ball, $|f(z)|\leq M|z|^k.$ 
\end{lemma}
Let $\mathcal{F}_{r,M,A}$ be the family of holomorphic mappings $f\colon r\mathbb{B}\rightarrow \mathbb{C}^N$, bounded by $M$, fixing the origin and with common differential $Az$ at the origin  satisfying $\|A\|<1$. 
\begin{lemma}\label{taylor}
For each $f\in \mathcal{F}_{r,M,A}$, we have  $|f(z)-Az|\leq C|z|^2,$ where $C=C(r,M,A).$
If moreover $f(z)-Az=O(|z|^k)$ for $k\geq 3$, then  $|f(z)-Az|\leq C_k|z|^k,$ where $C_k=C_k(r,M,A).$
\end{lemma}
\begin{proof}
Setting  $C_k=M/r^k+\|A\|/r^{k-1}$, the result follows from the previous lemma.
\end{proof}
As a consequence we get the following
\begin{lemma}\label{estimate}
For each $f\in \mathcal{F}_{r,M,A}$ we have the following estimate: to each $\|A\|<\alpha<1$ there corresponds $s>0$, $s=s(r,M,A)$ such that $|f(z)|\leq \alpha |z|,$ if $ |z|\leq s.$
\end{lemma}
\begin{proof} We proceed by contradiction: assume there exist a sequence $f_n\in\mathcal{F}_{r,M,A}$ and a sequence of points $z_n$ converging to the origin verifying $|f_n(z_n)|> \alpha |z_n|.$ We have $$|f_n(z_n)|=|Ax_n+f_n(z_n)-Az_n|\leq |Az_n|+ C|z_n|^2,$$ thus $$\alpha<\frac{|f_n(z_n)|}{|z_n|}\leq\frac{|Az_n|}{|z_n|}+C|z_n|,$$ but the right-hand term has $\limsup$ less or equal than $\|A\|$, which is the desired contradiction.
\end{proof}
\begin{lemma}\label{estimate2}
For each $f\in \mathcal{F}_{r,r,A}$ we have the following estimate: to each $s<r$ there corresponds  $K<1$, $K=K(r,A)$, such that $|f(z)|\leq K|z|,$ if $|z|\leq s.$ 
\end{lemma}
\begin{proof}
Assume the contrary: suppose there exist a sequence $f_n\in\mathcal{F}_{r,r,A}$ and a sequence of points $z_n$ in $\overline{s\mathbb{B}}$ verifying $|f_n(z_n)|> (1-1/n) |z_n|.$ Up to subsequences we have $z_n\rightarrow z'$ for some $z'$ such that $|z'|\leq s$, and $f_n\rightarrow f$ uniformly on compacta since $\mathcal{F}_{r,r,A}$ is a normal family. If $ z'\neq 0$ we have  $$1-\frac{1}{n}<\frac{|f_n(z_n)|}{|z_n|}\rightarrow \frac{|f(z')|}{|z'|},$$ and $|f(z')|/|z'|<1$ by Lemma \ref{schwarz}, which is a contradiction. If $ z'=0$, using again Lemma \ref{taylor} we get $$1-\frac{1}{n}<\frac{|f_n(z_n)|}{|z_n|}\leq\frac{|Az_n|}{|z_n|}+C|z_n|,$$ and the right-hand term has $\limsup$ less than or equal to $\|A\|$, contradiction.
\end{proof}
\begin{lemma}\label{koebe}
Suppose that $D$ is an open set in $\mathbb{C}^N$ containing the origin. Suppose we have an uniformly bounded family $\mathcal{H}$ of holomorphic mappings $h\colon D\rightarrow \mathbb{C}^N$ in $\Tang$. Then there exist a ball $r\mathbb{B}\subset D$ such that every $h\in\mathcal{H}$ is univalent on $r\mathbb{B}$, and a ball $s\mathbb{B}$ such that $s\mathbb{B}\subset h(r\mathbb{B})$ for all $h\in \mathcal{H}$.
\end{lemma}
\begin{proof}
Suppose there does not exist a ball  $r\mathbb{B}\subset D$ such that every $h\in\mathcal{H}$ is univalent on $r\mathbb{B}$. Since $\mathcal{H}$ is a normal family there exists a sequence $h_n\rightarrow  f$ uniformly on compacta, and such that  there does not exist a ball  $r\mathbb{B}\subset D$ with the property that every $h_n$ is univalent on $r\mathbb{B}$.
Since $f\in\Tang$  there exists a ball where $f$ is univalent. We can now apply \cite[Theorem 6.1.18]{Graham-Kohr}, getting a contradiction. Assume now there does not exist a ball contained in each $h(r\mathbb{B})$. Again there is a sequence $h'_n\rightarrow  f'$ uniformly on compacta, such that there does not exist a ball $s\mathbb{B}\subset\bigcap h_n'(r\mathbb{B}).$
The contradiction is then given by \cite[Proposition 3.1]{Arosio-Bracci-Hamada-Kohr}. 
\end{proof}

\section{Discrete evolution families and discrete Loewner chains}\label{discreteevolutionfamilies}
Let $A$ be a complex ($N\times N$)-matrix 
\begin{equation}
A=\mbox{Diag}(\lambda_1,\dots,\lambda_n),\ \ 0<|\lambda_N|\leq\cdots\leq |\lambda_1|<1.
\end{equation}
\begin{definition}
We define a {\sl discrete evolution family} on a domain $D\subseteq\mathbb{C}^N$ as a family $(\varphi_{n,m})_{ n\leq m\in \mathbb{N}}$ (sometimes denoted by $(\varphi_{n,m};D)$) of univalent self-mappings of $D$,  which satisfies the {\sl semigroup conditions}:  $$\varphi_{n,n}=\id,\ \ \ \varphi_{l,m}\circ\varphi_{n,l}=\varphi_{n,m},$$ where $0\leq n\leq l\leq m$. Such a family is clearly determined by the subfamily $(\varphi_{n,n+1})$. A {\sl dilation discrete evolution family} is a discrete evolution family such that
\begin{equation}\label{dilation type}
0\in D,\quad \varphi_{n,n+1}(0)=0,\quad \varphi_{n,n+1}(z)=Az+O(|z|^2)\ \mbox{for all}\  n\geq 0.
\end{equation}
\end{definition}
\begin{definition}
A family $(f_n)_{n\in\mathbb{N}}$ of holomorphic mappings $f_n\colon D\rightarrow \mathbb{C}^N$ is a {\sl discrete subordination chain} if for each $n<m$ the mapping $f_n$ is {\sl subordinate} to $f_m$, that is, there exists a holomorphic mapping (called {\sl transition mapping}) $\varphi_{n,m}\colon D\rightarrow D$ such that $$f_n=f_m\circ \varphi_{n,m}.$$ It is easy to see that the family of transition mappings of a subordination chain satisfies the semigroup property. If a subordination chain $(f_n)$ admits transition mappings  $\varphi_{n,m}$ which form a discrete evolution family (namely, $\v_{n,m}$ are univalent) we say that $(f_n)$ is {\sl associated} to $(\varphi_{n,m})$. 
\end{definition}
\begin{definition}
We define a {\sl discrete Loewner chain}  as a subordination chain $(f_n)$ such that every $f_n$ is univalent. In this case every {\sl transition mapping} $\varphi_{n,m}$ is univalent and uniquely determined. Thus the transition mappings form a discrete evolution family. A Loewner chain $(f_n)$ is {\sl normalized} if  $f_0(0)=0$ and $d_0f_0=\textrm{Id}$. A {\sl dilation discrete Loewner chain} is a discrete Loewner chain such that $$f_n(z)=A^{-n}z+O(|z|^2).$$ Following Pommerenke \cite{Pommerenke}, we call a dilation Loewner chain {\sl normal} if  $(A^n f_n)$ is a normal family. 
\end{definition}

\section{Triangular discrete evolution families}\label{triangulardiscrete}
Recall that a \textit{triangular automorphism} is a mapping $T\colon\mathbb{C}^N\rightarrow\mathbb{C}^N$ of the form 
$$
\begin{array}{l} 
T^{(1)}(z)=\lambda_1z_1, \\
T^{(2)}(z)=\lambda_2z_2+t^{(2)}(z_1), \\
T^{(3)}(z)=\lambda_3z_3+t^{(3)}(z_1,z_2), \\
\ \ \ \ \vdots \\
T^{(N)}(z)=\lambda_Nz_N+t^{(N)}(z_1,z_2,\dots,z_{N-1}),
\end{array} 
$$
where $0<|\lambda_j|<1$, and $t^{(i)}$ is a polynomial in $i-1$ variables, with all terms of degree greater or equal $2$. This is indeed an automorphism, since we can iteratively write its inverse, which is still a triangular automorphism:
\begin{equation}\label{triangularinverse}
\begin{array}{l} 
z_1=\frac{w_1}{\lambda_1}, \\
z_2=\frac{w_2}{\lambda_2}-\frac{1}{\lambda_2}t^{(2)}(z_1), \\
\ \ \ \ \vdots \\
z_N=\frac{w_N}{\lambda_N}-\frac{1}{\lambda_N}t^{(N)}(z_1,z_2,\dots,z_{N-1}).
\end{array} 
\end{equation}
\begin{definition}
The \textit{degree} of $T$ is $\max_i \deg T^{(i)}.$ We define a \textit{triangular evolution family} as a discrete dilation evolution family $(T_{n,m})$ of $\mathbb{C}^N$ such that each $T_{n,n+1}$, and hence every $T_{n,m}$, is a triangular automorphism. We denote $T^{(j)}_{n,n+1}(z)=\lambda_jz_j+t^{(j)}_{n,n+1}(z_1,z_2,\dots,z_{j-1})$ for all $n\geq 0$ and all $1\leq j\leq N$.  We denote $$T_{m,n}=T_{n,m}^{-1},\quad  0\leq n\leq m.$$ We say that a triangular evolution family $(T_{n,m})$ has {\sl bounded coefficients} if the family $(T_{n,n+1})$ has uniformly bounded coefficients, and we say that it has {\sl bounded degree} if $\sup_n \deg T_{n,n+1}<\infty$.
\end{definition} We can easily find a Loewner chain associated to a triangular evolution family: $$(f_n)=(T_{n,0})=(T_{0,n}^{-1}).$$ Indeed, $$f_m\circ T_{n,m}=T_{m,0}\circ T_{n,m}= T_{n,0}=f_n.$$   The following lemmas are just adaptations of  \cite[Lemma 1, p.80]{Rudin-Rosay}.

\begin{lemma}
Assume that $\sup_n \deg T_{n,n+1}<\infty$, then $\sup_n \deg T_{0,n}<\infty$.
\end{lemma}
\begin{proof}
Set for $j=1,\dots, N$, $$\mu^{(j)}=\max_n\deg T^{(j)}_{n,n+1}.$$ We denote by $S(m,k)$  the property $$\deg T^{(j)}_{0,k}\leq \mu^{(1)}\cdots\mu^{(j)},\ \ 1\leq j\leq m.$$ Since $T_{0,k+1}=T_{k,k+1}\circ T_{0,k},$ we have $$T^{(j)}_{0,k+1}=\lambda_jT^{(j)}_{0,k}+t^{(j)}_{k,k+1}(T^{(1)}_{0,k},\dots,T^{(j-1)}_{0,k}),\ \ 2\leq j\leq N,$$ thus $S(m,k+1)$ follows from  $S(m,k)$ and $S(m-1,k)$. Since  $S(1,k)$ and $S(m,1)$ are obviously true for all $k$ and $m$ (note that $\mu^{(1)}=1$, and $\mu^{(j)}\geq 1$, for every $j$),  $S(N,k)$ follows by induction.  Hence $$\deg T_{0,k}\leq \mu^{(1)}\cdots\mu^{(N)}.$$
\end{proof}

\begin{lemma}\label{rudinestimate}
Let $(T_{n,m})$ be a triangular evolution family of bounded degree and bounded coefficients. Let $\Delta$ be the unit polydisc. Then there exists a constant $\gamma>0$ such that $$T_{k,0}(\Delta)\subset \gamma^k\Delta,\ \ k\geq 1.$$ 
\end{lemma}
\begin{proof}
The family $(T_{n+1,n})$ of inverses of $(T_{n,n+1})$  has bounded coefficients. Indeed the family $(T_{n,n+1})$ has bounded coefficients, and the assertion follows by looking at (\ref{triangularinverse}). Likewise, $\sup_n \deg T_{n,0}<\infty$, since $\sup_n \deg T_{0,n}<\infty$.
Hence there exists $C\geq 1$ such that $|T_{n+1,n}^{(j)}(z)|\leq C$ for $z\in \Delta,\ 1\leq j\leq N,$ and there exists $d=\max_n \deg T_{n,0}.$  Let $M$ be the number of multi-indices $I=(i_1,\dots,i_N)$ with $|I|\leq d$, and set $\gamma=MC^d$, we claim that 
\begin{equation}\label{inductiveestimate}
|T^{(j)}_{k,0}(z)|\leq \gamma^k,\quad \mbox{for}\ z\in \Delta,\ j=1,\dots, N .
\end{equation} 
We proceed by induction on $k$.
Since $C\leq \gamma$, (\ref{inductiveestimate}) holds for $k=1$. Assume  (\ref{inductiveestimate}) holds for some $k\geq 1$. By Cauchy estimates the coefficients in $T^{(j)}_{k,0}(z)=\sum_{|I|\leq d}a_I z^I$ satisfy $$|a_{I}|\leq \gamma^k.$$ Since $T_{k+1,0}=T_{k,0}\circ T_{k+1,k},$ we have $$T^{(j)}_{k+1,0}=T^{(j)}_{k,0}(T^{(1)}_{k+1,k},\dots,T^{(N)}_{k+1,k})=\sum_{|I|\leq d}a_I(T^{(1)}_{k+1,k})^{i_1}\cdots (T^{(N)}_{k+1,k})^{i_N}.$$ Then $$|T^{(j)}_{k+1,0}|\leq MC^d\gamma^k=\gamma^{k+1}.$$ 
\end{proof}
\begin{corollary}\label{importante}
Let $(T_{n,m})$ be a triangular evolution family of  bounded degree and  bounded coefficients. Let $\frac{1}{2}\Delta$ be the polydisc of radius $(1/2)$.  Then there exists $\beta\geq 0$ such that for all $k\geq 1$ and all  $z,z'\in \frac{1}{2}\Delta$,  $$|T_{k,0}(z)-T_{k,0}(z')|\leq \beta^k|z-z'|.$$
\end{corollary}
\begin{proof}
Recall that $\Delta\subset \sqrt{N}\mathbb{B}$ and that if $B=(b_{ij})$ is a complex $(N\times N)$-matrix, then $$\|B\|\leq N\max_{i,j} |b_{ij}|.$$ If $z\in(1/2)\Delta$, then by Cauchy estimates and Lemma \ref{rudinestimate}, $$\|d_z T_{k,0}\|\leq 2N \sqrt N\gamma^k.$$

The result follows setting $\beta=  2N \sqrt N\gamma.$
\end{proof}
\begin{lemma}\label{uniqueness}
Let $(T_{n,m})$ be a triangular evolution family, with bounded degree and bounded coefficients. Then $T_{0,n}(z)\rightarrow 0$ uniformly on compacta. Hence for each neighborhood $V$ of $0$ we have $$\bigcup_{n=1}^\infty T_{n,0}(V)=\mathbb{C}^N.$$
\end{lemma}
\begin{proof}
Let $K$ be a compact set in $\mathbb{C}^N$. We proceed by induction on $i$. Notice that  $T^{(1)}_{0,k}(z)=\lambda_1^k z_1,$ hence if $\|\cdot\|$ denotes the sup-norm on $K$, we have $\|T^{(1)}_{0,k}\|\rightarrow 0$. Let $1< i\leq N$ and assume that $\lim_{k\rightarrow \infty}\|T^{(j)}_{0,k}\|=0,$ for $1\leq j<i.$ On $K$  $$\lim_{k\rightarrow \infty} \|t^{(i)}_{k,k+1}(T^{(1)}_{0,k},\dots,T^{(i-1)}_{0,k})\|=0,$$ since the family $(t^{(i)}_{k,k+1})$ has uniformly bounded coefficients and uniformly bounded degree. Notice that 
\begin{equation}\label{triangular}
T^{(i)}_{0,k+1}=\lambda_iT^{(i)}_{0,k}+t^{(i)}_{k,k+1}(T^{(1)}_{0,k},\dots,T^{(i-1)}_{0,k}),\ \ 2\leq i\leq N.
\end{equation}
Hence, for each  $\varepsilon >0$, $|T^{(i)}_{0,k+1}|\leq |\lambda_i\|T^{(i)}_{0,k}|+\varepsilon,$ on $K$ for $k$ large enough. Therefore 
$\limsup_{k\rightarrow \infty}\|T^{(i)}_{0,k}\|=0$, concluding the induction. 
\end{proof}
\section{Existence of discrete dilation Loewner chains:  nearly-triangular case}\label{existenceI}
We are going to prove the existence of Loewner chains associated to a given discrete dilation evolution family  by conjugating it to a triangular evolution family by means of a time dependent intertwining map. In this perspective, we shall see that normal Loewner chains correspond to time dependent linearizations of the evolution family.
\begin{definition}
Let $D\subset \mathbb{C}^N$ be a domain containing $0$.
Given two discrete dilation evolution families $(\varphi_{n,m};t\mathbb{B})$, $ (\psi_{n,m};D)$  suppose there exists, in a ball $r\mathbb{B}\subset t\mathbb{B}$, a normal family of univalent mappings $h_n\colon r\mathbb{B}\rightarrow D$  in $\Tang$, such that 
\begin{equation}\label{conjugation}
h_m\circ\varphi_{n,m}=\psi_{n,m}\circ h_n,\quad 0\leq n\leq m,
\end{equation}
then we shall say that $(h_n)$ {\sl conjugates} $(\varphi_{n,m})$ to $(\psi_{n,m})$.
\end{definition}
Notice that if $(\varphi_{n,m})$ is conjugate to $(\psi_{n,m})$ then necessarily $d_0\varphi_{n,m}=d_0\psi_{n,m}.$ 
Let $(\varphi_{n,m};t\mathbb{B})$ be a discrete dilation evolution family, and let $r\mathbb{B}\subset t\mathbb{B}$. Then $r\mathbb{B}$ is invariant for every $\varphi_{n,n+1}$, and hence $(\varphi_{n,m};t\mathbb{B})$ restricts to an evolution family $(\varphi_{n,m};r\mathbb{B})$.
\begin{lemma}\label{extension} 
Let $0<r<t$. Let $(\varphi_{n,m};t\mathbb{B})$ be a discrete dilation  evolution family. Suppose there exists a Loewner chain $(f_n)$ associated to the  evolution family $(\varphi_{n,m};r\mathbb{B}) $. Then there exists a  Loewner chain $(f_n^{\sf e})$ associated to $(\varphi_{n,m};t\mathbb{B})$  which extends $(f_n)$ in the following sense: $$f_n^{\sf e}(z)=f_n(z),\quad z\in r\mathbb{B},\ n\geq 0.$$
\end{lemma}
\begin{proof}
Fix $n\geq 0$. Let $0<r\leq s<t$. Let $k(s)\geq n$ be the least integer such that  $\v_{n,k(s)}(s\mathbb{B})\subset r\mathbb{B}$ (which exists by Lemma  \ref{estimate2}). Define $$f^{\sf e}_n(z)=f_{k(s)}(\v_{n,k(s)}(z)), \quad  z\in s\mathbb{B}.$$ A priori the value  $f^{\sf e}_n(z)$ depends on  $s$. However if  $0<r\leq s< u<t$,  then $k(u)\geq k(s)$, thus   $$f_{k(u)}(\varphi_{n,k(u)}(z))=f_{k(u)}(\varphi_{k(s),k(u)}(\varphi_{n,k(s)}(z)))=f_{k(s)}(\varphi_{n,k(s)}(z)), \quad \mbox{for all}\ z\in s\mathbb{B}.$$ Therefore $f^{\sf e}_n$ is well defined on $t\mathbb{B}.$
Notice that since $k(r)=n,$ $$f^{\sf e}_n|_{r\mathbb{B}}=f_n(\varphi_{n,n}(z))=f_n.$$  
It is easy to see that $f^{\sf e}_n$ is holomorphic and injective and that $(f^{\sf e}_n)$ is a Loewner chain associated to $(\varphi_{n,m};t\mathbb{B})$.
\end{proof}
Notice that the extended chain $(f^{\sf{e}}_n)$ can also be defined by 
\begin{equation}
f_n^{\sf e}(z)=\lim_{m\rightarrow\infty}f_m\circ\v_{n,m}(z),\quad z\in t\mathbb{B}.
\end{equation}
Now we can show how conjugations allow us to pull-back Loewner chains: 
\begin{remark}
Suppose that  $(h_n)$, defined on $r\mathbb{B}$,  conjugates $(\varphi_{n,m};t\mathbb{B})$ to $(\psi_{n,m};D)$, and assume that $(f_n)$ is a Loewner chain associated to $(\psi_{n,m})$. The {\sl pull-back chain} $(f_n\circ h_n)$ on $r\mathbb{B}$ is easily seen to be associated to $(\varphi_{n,m};r\mathbb{B})$.
By Lemma \ref{extension} one can extend $(f_n\circ h_n)$ to all of $t\mathbb{B}$ obtaining a Loewner chain  associated to $(\varphi_{n,m};t\mathbb{B})$ and defined by
\begin{equation}\label{pullbackextension}
  (f_n\circ h_n)^{\sf e}(z)=\lim_{m\rightarrow\infty}f_m\circ h_m\circ\v_{n,m}(z),\quad z\in t\mathbb{B}.
\end{equation}
\end{remark}

If $(h_n)$ conjugates an evolution family  $(\varphi_{n,m};t\mathbb{B})$  to a triangular evolution family $(T_{n,m})$, a Loewner chain associated to $(\varphi_{n,m})$ is given by the functions $(T_{n,0}\circ h_n)^{\sf e}$.
If in particular $(h_n)$ linearizes the given evolution family, that is  $T_{n,m}(z)=A^{m-n}z$, we obtain this way a normal Loewner chain $( A^{-n} h_n)$. Hence one has the following
\begin{proposition}\label{normaliff}
A discrete dilation evolution family $(\v_{n,m})$ admits a normal Loewner chain if and only if there exists a normal family $(h_n)$ of univalent mappings in $\Tang$ which conjugates it to its linear part: $$h_m\circ\varphi_{n,m}=A^{m-n} h_n, \quad 0\leq n\leq m.$$
\end{proposition}
Next we show how to find conjugations, provided we start with a discrete dilation evolution family  close enough to a triangular evolution family.
\begin{proposition}\label{koenigs}
Suppose that $(\varphi_{n,m};t\mathbb{B})$ is a discrete dilation evolution family, and that $(T_{n,m})$ is a triangular evolution family with bounded degree and bounded coefficients. Let $\beta$ be the constant given by Corollary \ref{importante} for $(T_{n,m})$, and let $k$ be an integer such that $$|\lambda_1|^k<\frac{1}{\beta}.$$  If for each $n\geq 0$ we have $$\varphi_{n,n+1}(z)-T_{n,n+1}(z)=O(|z|^k),$$ then $(\varphi_{n,m})$ is conjugate to $(T_{n,m})$.
\end{proposition}
\begin{proof}
Choose $|\lambda_1|<c<1$ such that $ c^k<1/\beta.$ Lemma \ref{estimate} gives us $r>0$ (we can assume $0<r<\min\{1/2,t\}$) such that on $r\mathbb{B}$ we have $|\varphi_{n,n+1}(z)|\leq c |z|$ and $|T_{n,n+1}(z)|\leq c |z|$ for all $n\geq 0$.   Thus for $z\in r\mathbb{B}$ we have $|\varphi_{0,n}(z)|<rc^n.$ Thanks to Lemma \ref{taylor} we have $$|\varphi_{n,n+1}(\zeta)-T_{n,n+1}(\zeta)|\leq C|\zeta|^k,$$ hence $$|\varphi_{0,n+1}(z)-T_{n,n+1}\circ\varphi_{0,n}(z)|=|\varphi_{n,n+1}\circ\varphi_{0,n}(z)-T_{n,n+1}\circ\varphi_{0,n}(z)| \leq C|\varphi_{0,n}(z)|^k\leq Cr^kc^{kn}.$$  The sequence $T_{n,0}\circ\varphi_{0,n}(z)$ verifies 
\begin{align}
|T_{n+1,0}\circ\varphi_{0,n+1}(z)-T_{n,0}\circ\varphi_{0,n}(z)| &= |T_{n+1,0}\circ\varphi_{0,n+1}(z)-(T_{n+1,0}\circ T_{n,n+1})\circ\varphi_{0,n}(z)| \nonumber \\
&\leq \beta^{n+1}   |\varphi_{0,n+1}(z)-T_{n,n+1}\circ\varphi_{0,n}(z)|   \nonumber \\
&\leq \beta^{n+1}   Cr^kc^{kn} \nonumber \\
&=\left(c^k\beta\right)^nCr^k\beta, \nonumber
\end{align}
where we used Corollary \ref{importante} (notice that since $r<1/2$, we have $|\varphi_{0,n+1}(z)|<c^{n+1}/2$ and $|T_{n,n+1}\circ\varphi_{0,n}(z)|<c^{n+1}/2,$ hence both $\varphi_{0,n+1}(z)$ and $T_{n,n+1}\circ\varphi_{0,n}(z)$ are  in $\frac{1}{2}\Delta$).

Hence the holomorphic mappings $T_{j,0}\circ\varphi_{0,j}(z)$  converge uniformly on $r\mathbb{B}$ to a holomorphic function $h_0\in\Tang$  (and univalent for the Hurwitz Theorem in several variables). Likewise, $$T_{j,n}\circ\varphi_{n,j}(z)\rightarrow h_n(z).$$ Each $h_n$ is bounded by $$1+\sum_{n=0}^\infty Cr^k\beta\left( c^k\beta\right)^n,$$ hence they form a normal family. Moreover $$h_m\circ\varphi_{n,m}=\lim_{j\rightarrow \infty} T_{j,m}\circ\varphi_{m,j}\circ\varphi_{n,m}=\lim_{j\rightarrow \infty}  T_{n,m}\circ T_{j,n}\circ\varphi_{n,j}= T_{n,m}\circ h_n.$$
\end{proof}

By (\ref{pullbackextension}), a Loewner chain on $t\mathbb{B}$ associated to $(\v_{n,m})$ is given by

\begin{equation}
(T_{n,0}\circ h_n)^{\sf e}(z)=\lim_{m\rightarrow\infty}T_{m,0}\circ(\lim_{j\rightarrow \infty} T_{j,m}\circ\varphi_{m,j})\circ \v_{n,m}(z)=\lim_{m\rightarrow\infty} T_{m,0}\circ \varphi_{n,m}(z).
\end{equation}

\section{Existence of discrete dilation Loewner chains:  general case}\label{existenceII}
In this section we show how to  conjugate a given discrete dilation evolution family  $(\varphi_{n,m};t\mathbb{B})$ to a nearly-triangular evolution family, by removing all non-resonant terms applying a parametric version of the Poincar\'e-Dulac method. This will give as a consequence the existence of Loewner chains for every discrete dilation evolution family.

\begin{definition}
A {\sl real resonance} for a matrix $A$ with eigenvalues $\lambda_i$ is an identity $$|\lambda_j|=|\lambda_1^{i_1}\dots\lambda_N^{i_N}|,$$ where $i_j\geq 0$, and $\sum_j i_j\geq 2$. If $|\lambda_j|<1$ for all $1\leq j\leq N$, real resonances can occur only in a finite number. Moreover, if $0<|\lambda_N|\leq\cdots\leq |\lambda_1|<1$ then the equality $|\lambda_j|=|\lambda_1^{i_1}\dots\lambda_N^{i_N}|,$ implies $i_j=i_{j+1}=\dots=i_N=0.$ 
Let $\v\colon t\mathbb{B}\rightarrow \mathbb{C}^N$ be a univalent mapping such that $\v(z)=Az+O(|z|^2)$, and denote its $j$-th component as  $$\varphi^{(j)}(z)=\lambda_j z_j+ \sum_{|I|\geq 2}a^{(j)}_{I} z^I,$$
where as usual, $z^I=z_1^{i_1}\cdots z_N^{i_N}$ for $I=(i_1,\dots, i_N).$ 
We call a monomial $a^{(j)}_{I} z^I$  {\sl resonant} if a real resonance $|\lambda_j|=|\lambda^I|$ occurs.  A mapping with only resonant monomials is necessarily triangular.
\end{definition}

\begin{proposition}\label{core}
Let $(\varphi_{n,m};t\mathbb{B})$ be a discrete dilation evolution family. For each $i\geq 2$  there exist
\begin{enumerate}
\item an evolution family  $(\varphi^{i}_{n,m})$ defined on a ball $\mathbb{B}_{i}$,
\item a uniformly bounded family $(k^{i}_n)$ of univalent maps defined on a ball $\mathbb{B}'_{i}\subset t\mathbb{B}$  which conjugates  $(\varphi_{n,m})$ to $(\varphi^{i}_{n,m})$,
\item a triangular evolution family $(T^i_{n,m})$  with $\deg T^i_{n,n+1}\leq i-1$ for all $n\geq 0$ and bounded coefficients such that for all $n\geq 0$,  $$\varphi^{i}_{n,n+1}=T^i_{n,n+1}+O(|z|^{i}).$$
\end{enumerate}
\end{proposition}
\begin{proof}
We proceed by induction. For $i=2$ it suffices to set $\varphi^2_{n,n+1}=\varphi_{n,n+1}$, $k^2_n=\id$ and $T^2_{n,n+1}=A$. Assume the proposition holds for $i\geq 2$.  Thus there exist  $(\varphi^{i}_{n,m} ;\mathbb{B}_{i})$   and $(T^i_{n,m})$ such that $$\varphi^{i}_{n,n+1}-T_{n,n+1}^i=O(|z|^i).$$ Since $\deg T^i_{n,n+1}\leq i-1$, we have $$\varphi^{i}_{n,n+1}-T_{n,n+1}^i-P^i_{n,n+1}=O(|z|^{i+1}),$$ where $P^i_{n,n+1}$ is the homogeneous term of $\varphi^{i}_{n,n+1}$ of degree $i$.  Let $R^i_{n,n+1}$ be the polynomial mapping obtained deleting every non-resonant term from $P^i_{n,n+1}$. Define the triangular evolution family $(T^{i+1}_{n,m})$ by $T^{i+1}_{n,n+1}=T^i_{n,n+1}+R^i_{n,n+1}$.  First we prove that there exists a family $(k_n)$ of   polynomial mappings in $\Tang$ with uniformly bounded degrees and uniformly bounded coefficients satisfying \begin{equation}\label{searching}
k_{n+1}\circ\varphi^{i}_{n,n+1}-T^{i+1}_{n,n+1}\circ k_n=O(|z|^{i+1}).
\end{equation} 
Let $I$ be a multi-index, $|I|=i$, and let $j$ be an integer $1\leq j\leq N$. Define $k_{I,j,n}$ as the polynomial mapping whose $l$-th component is $$k^{(l)}_{I,j,n}(z)= z_l +\delta_{lj}\alpha^{(j)}_{I,n}z^I,$$ where $\delta_{lj}$ is the Kronecker delta and $\alpha^{(j)}_{I,n}\in \mathbb{C}$ is to be chosen. Denote the $j$-th component of $\v_{n,n+1}$ as $\lambda_j z_j+ \sum_{|I|\geq 2}a^{(j)}_{I,n,n+1} z^I$. In the case $|\lambda_j|=|\lambda^I|$, that is when every $a^{(j)}_{I,n,n+1}z^I$ with $a^{(j)}_{I,n,n+1}\neq 0$ is resonant, set $\alpha^{(j)}_{I,n}=0$ for each $n$. In the case  $|\lambda_j|\neq |\lambda^I|$, by imposing the vanishing of terms in $z^I$ in the left-hand side of equation (\ref{searching}) we obtain the {\sl homological equation}:
\begin{equation}\label{homologicalequation}
\lambda^I \alpha^{(j)}_{I,n+1}+ a^{(j)}_{I,n,n+1}=\lambda_j \alpha^{(j)}_{I,n}.
\end{equation}
We have thus a recursive formula for  $\alpha^{(j)}_{I,n}$ in the non-resonant case: $$ \alpha^{(j)}_{I,n+1}=\frac{\lambda_j \alpha^{(j)}_{I,n}-a^{(j)}_{I,n,n+1}}{\lambda^I},$$ which gives
\begin{equation}\label{magic}\alpha^{(j)}_{I,n}=\alpha^{(j)}_{I,0}\left(\frac{\lambda_j}{\lambda^I}\right)^n-\frac{a^{(j)}_{I,0,1}}{\lambda^I}\left(\frac{\lambda_j}{\lambda^I}\right)^{n-1}-\frac{a^{(j)}_{I,1,2}}{\lambda^I}\left(\frac{\lambda_j}{\lambda^I}\right)^{n-2}-\cdots-\frac{a^{(j)}_{I,n-1,n}}{\lambda^I}.
\end{equation}  
Since by Cauchy estimates $(a^{(j)}_{I,n,n+1})$ is a bounded sequence, if $|\lambda_j|<|\lambda^I|$ then $(\alpha^{(j)}_{I,n})$ is bounded regardless of our choice for $\alpha^{(j)}_{I,0}\in \mathbb{C}$, so we can set $\alpha^{(j)}_{I,0}=0$. In the case  $|\lambda_j|>|\lambda^I|$ we have to choose $\alpha^{(j)}_{I,0}$ suitably in order to obtain a bounded sequence.
Divide (\ref{magic}) by $(\lambda_j/\lambda^I)^n$:
\begin{equation}\label{magic2}
\alpha^{(j)}_{I,n}\left(\frac{\lambda^I}{\lambda_j}\right)^n=\alpha^{(j)}_{I,0}-\frac{a^{(j)}_{I,0,1}}{\lambda^I}\left(\frac{\lambda^I}{\lambda_j}\right)-\frac{a^{(j)}_{I,1,2}}{\lambda^I}\left(\frac{\lambda^I}{\lambda_j}\right)^2-\cdots-\frac{a^{(j)}_{I,n-1,n}}{\lambda^I}\left(\frac{\lambda^I}{\lambda_j}\right)^n,
\end{equation} 
and set $$\alpha^{(j)}_{I,0}=\sum_{m=1}^\infty \frac{a^{(j)}_{I,m-1,m}}{\lambda^I}\left(\frac{\lambda^I}{\lambda_j}\right)^m,$$ which converges since  $(a^{(j)}_{I,n,n+1})$ is a bounded sequence.
With this choice, $$\alpha^{(j)}_{I,n}\left(\frac{\lambda^I}{\lambda_j}\right)^n=\sum_{m=n+1}^\infty \frac{a^{(j)}_{I,m-1,m}}{\lambda^I}\left(\frac{\lambda^I}{\lambda_j}\right)^m,$$ so that $$ \alpha^{(j)}_{I,n}=\sum_{m=n+1}^\infty \frac{a^{(j)}_{I,m-1,m}}{\lambda^I}\left(\frac{\lambda^I}{\lambda_j}\right)^{m-n},$$ hence $(\alpha^{(j)}_{I,n})$ is also bounded. Fix an order on the set $\{(I,j):|I|=i,\ 1\leq j\leq N\}$ and define the mapping $k_n$ as the ordered composition of all $k_{I,j,n}$ with $|I|=i,\ 1\leq j\leq N$. It is easy to check that $(k_n)$ is a family of  polynomial mappings in $\Tang$ with uniformly bounded degree and uniformly bounded coefficients  satisfying  (\ref{searching}).

Lemma \ref{koebe} yields a ball $r\mathbb{B}\subset D$ such that every $k_n$ is univalent on $r\mathbb{B}$, and a ball $s\mathbb{B}$ such that $s\mathbb{B}\subset k_n(r\mathbb{B})$ for all $n\geq 0$. On $s\mathbb{B}$ we can define a family of holomorphic mappings $$\varphi^{i+1}_{n,n+1}=k_{n+1}\circ \varphi^{i}_{n,n+1}\circ k_n^{-1}.$$ By Lemma \ref{estimate} there exists a ball $\mathbb{B}_{i+1}$ invariant for each $\varphi^{i+1}_{n,n+1}$. Hence  $(\varphi^{i+1}_{n,n+1};\mathbb{B}_{i+1})$ is a discrete evolution family. Since $(k_n)$ is an equicontinuous family, there exists a ball  $u\mathbb{B}$ such that  $k_n(u\mathbb{B})\subset \mathbb{B}_{i+1}$ for all $n\geq 0$, so that $(k_n)$ conjugates $(\varphi^{i}_{n,n+1};\mathbb{B}_{i})$ to $(\varphi^{i+1}_{n,n+1};\mathbb{B}_{i+1})$:  $$k_{n+1}\circ \varphi^{i}_{n,n+1}=\varphi^{i+1}_{n,n+1}\circ k_n.$$ Since (\ref{searching}) holds by construction, we have $$\varphi^{i+1}_{n,n+1}\circ k_n-T^{i+1}_{n,n+1}\circ k_n=O(|z|^{i+1}),$$ that is, $$\varphi^{i+1}_{n,n+1}-T^{i+1}_{n,n+1}=O(|z|^{i+1}).$$ The family $(k_n^{i+1})$ conjugating $(\varphi_{n,m})$ to $(\varphi^{i+1}_{n,m})$ is obtained by composing for each $n\geq 0$ the conjugation mappings $k_n^{i}$ given by inductive hypothesis with the conjugation mappings  $k_n$. Indeed since the family $(k_n^{i})$ is equicontinuous by inductive hypothesis, there exists a ball $\mathbb{B}'_{i+1}$ such that $k^{i}_n(\mathbb{B}'_{i+1})\subset u\mathbb{B}$ for all $n\geq 0$. Let  $(k^{i+1}_n)$ be the family of mappings defined on $\mathbb{B}'_{i+1}$ by $(k_n\circ k^{i}_n)$, then

\begin{align}
k^{i+1}_{n+1}\circ \varphi_{n,n+1}=(k_{n+1}\circ k^{i}_{n+1})\circ \varphi_{n,n+1}&=k_{n+1}\circ\varphi^{i}_{n,n+1}\circ k^{i}_n\nonumber \\
&=\varphi^{i+1}_{n,n+1} \circ(k_n\circ k^{i}_n)\nonumber \\
&=\varphi^{i+1}_{n,n+1}\circ k^{i+1}_{n}.\nonumber
\end{align}
 \end{proof}

\begin{remark}
Let $q$ be the smallest integer such that $|\lambda_1^q|<|\lambda_N|.$ Then no term of $P^q_{n,n+1}$ can be resonant. Hence $T^i_{n,n+1}=T^q_{n,n+1}$ for any $i\geq q$.
\end{remark}

\begin{proposition}\label{result}
A discrete dilation  evolution family $(\varphi_{n,m};t\mathbb{B})$ admits an associated normalized Loewner chain $(f_n)$ such that $\bigcup_n f_n(t\mathbb{B})=\mathbb{C}^N.$ If no real resonances occur $(f_n)$ is a normal chain.
\end{proposition}

\begin{proof} 
Denote $(T_{n,m})=(T^q_{n,m}),$ where $q$ is the smallest integer such that $|\lambda_1^q|<|\lambda_N|.$
Let $\beta$ be the constant given by Lemma \ref{inductiveestimate} for $(T_{n,m})$, and let $l$ be an integer such that $|\lambda_1|^l<\frac{1}{\beta}.$ Let $(\varphi^{l}_{n,m};\mathbb{B}_{l})$ be the evolution family given by Proposition \ref{core}. We can apply Proposition \ref{koenigs} obtaining  a uniformly bounded family $(h_n)$ given by $$h_n=\lim_{m\rightarrow \infty}T_{m,n}\circ\varphi^l_{n,m},$$ defined on a ball $r\mathbb{B}\subset \mathbb{B}_{l}$, which conjugates  $(\varphi^{l}_{n,m};\mathbb{B}_{l})$ to $(T_{n,m})$. 

Thus a Loewner chain associated to $(\varphi^{l}_{n,m};\mathbb{B}_{l})$ is given by the mappings $$(T_{n,0}\circ h_n)^{\sf e}=\lim_{m\rightarrow \infty} T_{m,0}\circ\varphi^{l}_{n,m}.$$  Since $(k_n^l)$ conjugates  $(\varphi_{n,m};t\mathbb{B})$ to $(\varphi^{l}_{n,n+1};\mathbb{B}_l)$, a Loewner chain associated to  $(\varphi_{n,m};t\mathbb{B})$ is given by  $$f_n=((T_{n,0}\circ h_n)^{\sf e}\circ k_n^{l})^{\sf e}=\lim_{m\rightarrow \infty} T_{m,0}\circ k^{l}_{m}\circ \varphi_{n,m}.$$ 

Now we prove $\bigcup_n f_n(t\mathbb{B})=\mathbb{C}^N$. Since the family $(k_n^l)$ is equicontinuous, there exists a ball $u\mathbb{B}\subset t\mathbb{B}$ such that $k^l_n(u\mathbb{B})\subset r\mathbb{B}$ for all $n\geq 0$. On $u\mathbb{B}$,  $$f_n= T_{n,0}\circ h_n\circ k_n^l.$$
The sequence $h_n \circ k^l_n$ is uniformly bounded, hence  Lemma \ref{koebe} yields the existence of  a ball $V$ contained in each $h_n\circ k^l_n(u\mathbb{B})$. Thus $$\bigcup_n f_n(t\mathbb{B})\supseteq \bigcup_n T_{n,0}(V)=\mathbb{C}^N.$$

If no real resonances occur then $T_{n,m}(z)=A^{m-n}z,$ hence $$f_n=\lim_{m\rightarrow \infty} A^{-m} k^l_{m}\circ \varphi_{n,m}.$$  As above we have that on $u\mathbb{B}$, the sequence $$A^n f_n= \lim_{m\rightarrow \infty}\left( A^{n-m}\varphi^l_{n,m}\right) \circ k^l_n$$ is uniformly bounded. Let $s\mathbb{B}$ be a ball contained in $t\mathbb{B}$. Lemma \ref{estimate2} yields an integer $j_n$ such that  $\varphi_{n,j_n}(s\mathbb{B})\subset u\mathbb{B}$ and $j_n-n$ does not depend on $n$. From $$A^n f_n=A^{n-j_n} A^{j_n} f_{j_n}\circ\v_{n,j_n}$$ we see that $A^n f_n$ is uniformly bounded on $s\mathbb{B}$, hence it is a normal family.
\end{proof}

\section{Essential uniqueness}\label{essentialuniqueness}
\begin{proposition}\label{uniquenessresult}
Let $(\varphi_{n,m};t\mathbb{B})$ be a discrete dilation evolution family, and let $(f_n)$ be the Loewner chain given by Proposition \ref{result}.  A family of holomorphic mappings $(g_n)$ is a subordination chain  associated to $(\varphi_{n,m})$ if and only if there exists an entire mapping $\Psi$ on $\mathbb{C}^{N}$ such that $$g_n=\Psi\circ f_n.$$ 
\end{proposition}
\begin{proof}
Set $\Psi_n= g_n \circ f_n^{-1}.$ If $m>n$,
\begin{equation}\label{azz}
\Psi_m|_{f_n(t\mathbb{B})}=\Psi_n,
\end{equation}
as it is clear from the following commuting diagram
$$\xymatrix{\mathbb{C}^N\ar[r]^{\id}& \mathbb{C}^N\\
t\mathbb{B}\ar[r]^{\varphi_{n,m}}\ar[u]^{f_n}\ar[d]_{g_n}& t\mathbb{B}\ar[u]^{f_m}\ar[d]_{g_m}\\
\mathbb{C}^N\ar[r]^{\id}& \mathbb{C}^N.}$$
Therefore  by (\ref{azz}) we can define on $\mathbb{C}^N=\bigcup_n f_n(t\mathbb{B})$ a  mapping $\Psi$ setting $$ \Psi|_{f_n(t\mathbb{B})}=\Psi_n.$$ This proves the first statement, and the converse is trivial.
\end{proof}

\section{Continuous case}\label{continuouscase}

Let $\Lambda $ be a complex ($N\times N$)-matrix 
\begin{equation}
\Lambda =\mbox{Diag}(\alpha_1,\dots,\alpha_N),\quad \mbox{where}\quad\Re{\alpha_N}\leq\dots\leq\Re{\alpha_1}<0.
\end{equation}
\begin{definition}
We define a {\sl dilation evolution family} as a family $(\v_{s,t})_{ 0\leq s\leq t}$ of holomorphic self-mappings of the unit ball $\mathbb{B}\subset \mathbb{C}^N$ such that for any $0\leq s\leq u \leq t$, $$\v_{s,s}=\id_{\mathbb{B}},\quad\v_{s,t}=\v_{u,t}\circ\v_{s,u},\quad \varphi_{s,t}(z)=e^{\Lambda ( t-s)}z+O(|z|^2).$$
\end{definition}
\begin{definition}
A family $(f_s)_{s\geq 0}$ of holomorphic mappings $f_s\colon\mathbb{B}\rightarrow \mathbb{C}^N$ is called a {\sl subordination chain} if $f_s$ is subordinated to $f_t$ when $s\leq t$. If a subordination chain admits as transition mappings a dilation evolution family $(\varphi_{s,t})$ we say that $(f_s)$ is {\sl associated} to $(\varphi_{s,t})$. We define a {\sl dilation Loewner chain} as a subordination chain such that each $f_s$ is univalent and $$f_s(z)=e^{-\Lambda s}z+O(|z|^2).$$ A dilation Loewner chain is {\sl normal} if $(e^{\Lambda s}f_s)$ is a normal family.
\end{definition}
\begin{definition}
 We define a {\sl dilation Herglotz vector field} $H(z,t)$ as a function $H\colon\mathbb{B}\times [0,+\infty)\rightarrow \mathbb{C}^N$ such that  for all $z \in\mathbb{B}$ the mapping $H(z,\cdot)$ is measurable, and such that $H(\cdot,t)$ is a holomorphic mapping  for  a.e. $t\geq 0$ satisfying
$$H(z,t)=\Lambda z + O(|z|^2),\quad \Re\langle H(z,t),z\rangle \leq 0\quad \forall z\in\mathbb{B}.$$
\end{definition}
\begin{lemma}\label{lemmadilation}
Let $k_\mathbb{B}$ be the Kobayashi metric of $\mathbb{B}.$
Given a dilation evolution family $(\v_{s,t})$, a dilation Loewner chain $(f_s)$ and a dilation Herglotz vector field $H(z,t)$ the following hold: to each $T>0$ and to any compact set $K\subset \mathbb{B}$ there correspond positive constants $c_{T,K},C_{T,K}$ and $k_{T,K}$ such that for all $z\in K$ and $0\leq s\leq t' \leq t\leq T,$
\begin{enumerate}
\item $k_{\mathbb{B}}(\varphi_{s,t}(z),\varphi_{s,t'}(z))\leq c_{T,K}(t-t'),$
\item $|f_s(z)-f_t(z)|\leq k_{K,T}(t-s),$
\item $|H(z,t)|\leq C_{K,T}.$
\end{enumerate}
Therefore $(\v_{s,t})$ is an $L^\infty$-evolution family, $(f_s)$ is an $L^\infty$-Loewner chain, and $H(z,t)$ is an $L^\infty$-Herglotz vector field, in the sense of \cite{Arosio-Bracci-Hamada-Kohr}\cite{Bracci-Contreras-Diaz-II}.

\end{lemma}
\begin{proof}
Let $H(z,t)$ be a dilation Herglotz vector field. By Lemma 1.2 \cite{Graham-Hamada-Kohr-Kohr}, we have on $r\mathbb{B}$ $$|H(z,t)|\leq \frac{4r}{(1-r)^2}\|\Lambda \|.$$ Hence $H(z,t)$ is an $L^\infty$-Herglotz vector field. For  $(\varphi_{s,t})$ and $(f_s)$, see the proof of \cite[Theorem 2.8]{Graham-Hamada-Kohr-Kohr}.
\end{proof}
Recall if $(\v_{s,t})$ is an  $L^\infty$-evolution family, each mapping $\v_{s,t}$ is univalent   \cite[Proposition 5.1]{Bracci-Contreras-Diaz-II}. If we restrict the time to integer values in a dilation evolution family $(\varphi_{s,t})$  we obtain the {\sl discretized} dilation evolution family $(\varphi_{n,m})$. We have $Az=d_0\varphi_{n,n+1}(z)= (e^{\alpha_1}z_1,\dots,e^{\alpha_N}z_N)=e^{\Lambda }z.$ In the continuous framework a {\sl real resonance} is an identity $$\Re(\sum_{j=1}^N k_j\alpha_j)=\Re\alpha_l,$$ where $k_j\geq0$ and $\sum_j k_j\geq 2$. It is easy to see that a continuous real resonance corresponds to a real resonance for the discretized evolution family.

\begin{lemma}\label{discretetocontinuous}
Let $(\varphi_{s,t})$ be a dilation evolution family, and let $(\varphi_{n,m})$ be its discretized evolution family. Assume there exists a discrete Loewner chain $(f_n)$ associated to  $(\varphi_{n,m})$. Then we can extend it in a unique way to a dilation Loewner chain associated to $(\varphi_{s,t})$. If $(f_n)$ is a normal Loewner chain, then also $(f_s)$ is  normal.
\end{lemma}
\begin{proof}
Define for $s\geq 0$, $$f_s=f_j\circ \varphi_{s,j},$$ where $j$ is an integer such that $s\leq j$. The mapping $f_s$ is well defined. Indeed, let $j<k$, then $$f_j\circ\varphi_{s,j}=f_k\circ \varphi_{j,k}\circ\varphi_{s,j}=f_k\circ\varphi_{s,k}.$$ The family $(f_s)$ is a subordination chain: indeed if $0\leq s\leq t\leq j$, then $$f_s=f_j\circ\varphi_{s,j}=f_j\circ\varphi_{t,j}\circ\varphi_{s,t}=f_t\circ \varphi_{s,t}.$$ Moreover each $f_s$ is univalent and $d_0f_s(z)=e^{-\Lambda s}z$, thus Lemma \ref{lemmadilation} yields that $(f_s)$ is a dilation $L^\infty$-Loewner chain. Assume  $(e^{\Lambda n} f_n)$ is a normal family. If $0<r<1$ this family is uniformly bounded on $r\mathbb{B}$. For each $s\geq 0$ define $m_s$ as the smallest integer greater than $s$. We have $$ e^{\Lambda s}f_s=e^{\Lambda s} f_{m_s}\circ\v_{s,m_s}= e^{\Lambda (s-m_s)} e^{\Lambda m_s} f_{m_s} \circ\v_{s,m_s},$$
which is uniformly bounded on $r\mathbb{B}$ since $m_s-s$ is smaller than $1$. Hence $(e^{\Lambda s} f_s)$ is a normal family.
\end{proof}
The following result generalizes Theorem 2.3 in \cite{Graham-Hamada-Kohr-Kohr} and Theorem 3.1 in \cite{Duren-Graham-Hamada-Kohr}, where the hypothesis $2 \Re \alpha_1< \Re \alpha_N$ implies that no real resonances can occur (however in such papers the authors consider  non-necessarily diagonal $\Lambda$).
\begin{theorem}\label{continuousresult}
Let $(\varphi_{s,t})$ be a dilation evolution family. Then there exists a dilation Loewner chain $(f_s)$ associated to $(\varphi_{s,t})$, such that $\bigcup_s f_s(\mathbb{B})=\mathbb{C}^N.$ If no real resonances occur then $(f_s)$ is a normal chain.  A family of holomorphic mappings $(g_s)$ is a subordination chain  associated to $(\varphi_{s,t})$ if and only if there exists an entire mapping $\Psi$ on $\mathbb{C}^{N}$ such that $$g_s=\Psi\circ f_s.$$ 
\end{theorem}
\begin{proof}
The result follows by applying Proposition \ref{result} to the discretized evolution family associated to $(\varphi_{s,t})$, then Lemma \ref{discretetocontinuous} and Proposition \ref{uniquenessresult}.
\end{proof}
\begin{remark}
For $(f_s)$ we have the expression (with notations as in Proposition \ref{result}) $$f_s(z)= \lim_{m\rightarrow \infty} T_{m,0}\circ k^l_{m}\circ \varphi_{s,m}(z).$$
If we assume  $2 \Re \alpha_1<\Re \alpha_N$ then no real resonances can occur. Thus in this case $T_{m,0}(z)=A^{-m}z=e^{-\Lambda m}z$, so that the constant $\beta$ given by Corollary \ref{importante} can be taken equal to $\|A^{-1}\|=1/\lambda_N$. Hence $|\lambda_1|^l<1/\beta$ holds for $l=2$.  Therefore we can use Proposition \ref{koenigs} directly, obtaining $$f_s(z)= \lim_{m\rightarrow \infty}e^{-\Lambda m} \varphi_{s,m}(z),$$ in agreement with  Theorem 2.3  of \cite{Graham-Hamada-Kohr-Kohr}.
\end{remark}
Recall \cite[Theorem 2.1]{Graham-Hamada-Kohr-Kohr} that if $H$ is a dilation  Herglotz vector field, the {\sl Loewner ODE}
\begin{equation}
\begin{cases}
\overset{\bullet}z(t)=H(z,t)\\
z(s)=z
\end{cases}
\end{equation}
has a unique solution $t\mapsto\v_{s,t}(z)$, and $(\v_{s,t})$ is a dilation $L^\infty$-evolution family.
\begin{definition}
The partial differential equation $$\frac{\partial f_t(z)}{\partial t}=-d_zf_tH(z,t)\quad t\geq 0,\ z\in\mathbb{B}, $$ where $H(z,t)$ is a dilation  Herglotz vector field, is called the {\sl Loewner PDE}.
\end{definition}
With these notations, Theorem \ref{continuousresult} can be rephrased as 
\begin{theorem}
Let $H$ be a dilation Herglotz vector field, and let $t\mapsto\v_{s,t}$ be the solution of the associated Loewner ODE. Then if $(f_s)$ is the Loewner chain associated to the dilation $L^\infty$-evolution family $(\v_{s,t})$ given by Theorem \ref{continuousresult}, the mapping $t\mapsto f_t$ is a solution for the  Loewner PDE $$\frac{\partial f_t(z)}{\partial t}=-d_zf_tH(z,t).$$ Moreover,  a family $(g_s)$ of holomorphic mappings on the ball satisfies
\begin{enumerate}
\item the mapping $t\mapsto g_t$ is locally absolutely continuous in $t$, uniformly on compacta with respect to $z\in \mathbb{B}$,
\item  the mapping $t\mapsto g_t$ solves the Loewner PDE,
\end{enumerate}
if and only if there exists an entire mapping $\Psi$ on $\mathbb{C}^{N}$ such that $$g_s=\Psi\circ f_s.$$ 

\end{theorem}
\begin{proof}
It suffices to recall that such a $t\mapsto g_t$ satisfies the Loewner PDE if and only if $(g_s)$ is a subordination chain associated to  $(\v_{s,t})$. See the proof of \cite[Theorem 3.1]{Duren-Graham-Hamada-Kohr}.
\end{proof}

\begin{remark}
A dilation evolution family $(\v_{s,t})$ is called {\sl periodic} if $\v_{s,t}=\v_{s+1,t+1},$  for all $0\leq s\leq t$.
For periodic dilation evolution families the {\sl pure real resonances}, that is real resonances which are not  complex resonances,   are not obstructions to the existence of a normal Loewner chain. Namely if $(\v_{s,t})$ is a periodic dilation evolution family and no complex resonances occur, then there exists a  normal Loewner chain $(f_s)$ associated to $(\v_{s,t}).$ Indeed it is easy to see that by the Poincar\'e Theorem \cite[pp. 80--86]{Rudin-Rosay} the discretized evolution family $(\v_{n,m})=(\v_{0,1}^{\circ(m-n)})$ admits a discrete normal Loewner chain $(f_n)$. Lemma \ref{discretetocontinuous} yields then a normal Loewner chain $(f_s).$
\end{remark}

\section{Counterexamples}\label{counterexamples}

\subsection*{1}\label{counter1}
Let $\Lambda =\mbox{Diag}(\alpha_1,\alpha_2)$. If $(\varphi_{s,t})$ is a dilation evolution family such that $\v_{s,t}=e^{\Lambda(t-s)}z+O(|z|^2)$ and $2 \Re \alpha_1< \Re \alpha_N$, then by Lemma 2.12 in \cite{Duren-Graham-Hamada-Kohr} there exists a unique normal Loewner chain associated to  $(\varphi_{s,t})$. This is no longer true when $2 \Re \alpha_1\geq \Re \alpha_N$. Indeed, consider on $\mathbb{B}\subset \mathbb{C}^2$ the linear dilation evolution family  defined by $$\v_{s,t}(z)=e^{\Lambda (t-s)}z=(e^{\alpha_1(t-s)}z_1,e^{\alpha_2(t-s)}z_2).$$ The family $(e^{-\Lambda s}z)$ is trivially a normal Loewner chain associated to $(e^{\Lambda (t-s)}z)$. The univalent family $$k_s(z)=(z_1,z_2+ e^{(\alpha_2-2\alpha_1)s}z_1^2),$$ satisfies $k_{t}\circ e^{\Lambda (t-s)}z=e^{\Lambda (t-s)} k_{s}(z).$ Since $\Re\alpha_2\leq 2\Re \alpha_1$, it is a uniformly bounded family, thus $(e^{-\Lambda s} k_s)$ is another normal Loewner chain associated to $(e^{\Lambda (t-s)}z)$.

\subsection*{2}\label{counter2}
Let  $\Lambda =\mbox{Diag}(\alpha_1,\alpha_2)$, $\alpha_2=2\alpha_1$. There exists a dilation evolution family $(\v_{s,t})$ such that $\v_{s,t}=e^{\Lambda(t-s)}z+O(|z|^2)$, which does not admit any associated normal Loewner chain.
Indeed,  for $c\in\mathbb{C}^*$ small enough, the family $(\psi_t)$ defined by $$\psi_t(z)=(e^{\alpha_1 t}z_1,e^{\alpha_2t}(z_2+ctz_1^2))$$ is a semigroup on $\mathbb{B}\subset\mathbb{C}^2$. Thus  $$\varphi_{s,t}(z)=\psi_{t-s}(z)$$ defines a dilation evolution family.  Assume by contradiction there exists a normal Loewner chain $(f_s)$  associated to $(\v_{s,t})$. The family $(h_s)=(e^{\Lambda s}f_s)$ satisfies $h_t\circ\varphi_{s,t}=e^{\Lambda (t-s)} h_s,$ so in particular  
\begin{equation}\label{enlc}
h_t\circ\varphi_{0,t}=e^{\Lambda t} h_0.
\end{equation}
Let $a_s$ be the coefficient of the term $z_1^2$ in the second component of $h_s$. Then imposing equality of terms in $z_1^2$ in equation (\ref{enlc}) we find $e^{\alpha_2t}ct+a_te^{2\alpha_1t}= a_0e^{\alpha_2t},$ hence $$a_t=e^{(\alpha_2-2\alpha_1)t}(a_0-ct),$$ which gives
$a_t=a_0-ct,$ so that $(h_s)$ cannot be a normal family.

\subsection*{3}\label{counter2b}
Let $A=\mbox{Diag}(\lambda_1,\lambda_2)$,  $|\lambda_1|^2=|\lambda_2|$, $\lambda_1^2\neq \lambda_2$.  There exists a discrete dilation evolution family $(\v_{n,m})$ such that $\v_{n,n+1}(z)=Az+O(|z|^2)$, which does not admit any associated discrete normal Loewner chain. Indeed, if $r>0$ is sufficiently small, given any sequence $(a_{n,n+1})$ in $r\mathbb{B}$ there exists a discrete dilation evolution family defined by $$\v_{n,n+1}(z)=(\lambda_1 z_1, \lambda_2 z_2 +a_{n,n+1}z_1^2).$$  Assume by contradiction there exists a normal Loewner chain $(f_n)$  associated to $(\v_{n,m})$. The family $(h_n)=(A^n f_n)$ satisfies 
\begin{equation}\label{buh}
h_{n+1}\circ\varphi_{n,n+1}=A h_n.
\end{equation}
Let $\alpha_n$ be the coefficient of the term $z_1^2$ in the second component of $h_s$, and set $\zeta=\lambda_1^2/\lambda_2$. Then imposing equality of terms in $z_1^2$ in equation (\ref{buh}) we obtain as in (\ref{magic2}) the recursive formula 
$$\alpha_n\zeta^n\lambda_1^2=\alpha_0\lambda_1^2-a_{0,1}\zeta-a_{1,2}\zeta^2-\cdots-a_{n-1,n}\zeta^n.$$ 
For $1\leq j\leq 8$ define $C_j=\{\zeta\in S^1 : 2\pi(j-1)/8\leq\arg z\leq 2\pi j/8\} $. There exists a $C_j$ which contains the images of a subsequence $(\zeta^{k_n})$. Set 
\begin{equation}
a_{m-1,m}=
\begin{cases}
r/2,\  \mbox{if there exists $n$ such that $m=k_n$},\\ 
0,\  \mbox{otherwise},
\end{cases}
\end{equation}
then the sequence $(\sum_{j=0}^n a_{j-1,j}\zeta^j)$ is not bounded, hence the sequence $(\alpha_n)$ is also not  bounded.
Thus for $(\v_{n,m})$ no normal family $(h_n)$ can solve (\ref{buh}).

\subsection*{4}\label{counter3}
There exists  a discrete evolution family $(\v_{n,m})$ on $\mathbb{B}^3\subseteq\mathbb{C}^3$ which does not admit any associated discrete Loewner chain. Indeed, by \cite{Fornaess} there exists a complex manifold $M$ which is an increasing union of open sets $M_n$ each of which biholomorphic to the ball $\mathbb{B}^3$ by means of a biholomorphism $f_n\colon \mathbb{B}^3 \to M_n$, with the property that $M$ is not Stein. By  \cite{Behnke-Stein} this implies that $M$ cannot be embedded into $\mathbb{C}^3$ as an open set.

Define $\v_{n,n+1}=  f_{n+1}^{-1}\circ f_n$ for all $n\geq 0$. Then $(\v_{n,m})$ is a discrete evolution family which does not admit any associated discrete Loewner chain $(g_n)$. Indeed, if such a family existed, then $M$ would be biholomorphic to the open subset of $\mathbb{C}^3$ given by $\bigcup_n g_n(\mathbb{B}^3).$

This suggests the following (open) question: does such a discrete evolution family embed into some $L^\infty$-evolution family on $\mathbb{B}^3$?

\end{document}